\documentclass{amsart}
\usepackage{graphicx}
\usepackage{tikz-cd}
\usetikzlibrary{cd}
\newtheorem{theorem}{Theorem}

\newcommand{\abs}[1]{\left\vert#1\right\vert}

\usepackage{amsmath,amscd}
\begin{document}
\title{Explicit Holomorphic Structures for Embeddings of Closed 3-manifolds into $\mathbb{C}^3$}
\author{Ali M. Elgindi
\\ Institute of Mathematics, Henan Academy of Sciences. 
\\E-mail: ali.m.elgindi@gmail.com}

% ----------------------------------------------------------------

\begin{abstract}  
Expanding on my former work along with the more recent work of Kasuya and Takase, we demonstrate that for a given link $L \subset M$ which is null-homologous in $H_1(M)$ and for any smooth oriented 2-plane field $\eta$ over $L$ there exists a smooth embedding $F:M \hookrightarrow \mathbb{C}^3$ so that the set of complex tangents to the embedding is exactly $L$ and at each $x \in L$ the holomorphic tangent space is exactly $\eta_x$. Furthermore, we demonstrate how the analyticity of a complex tangent, as given by the Bishop invariant, may be determined exactly from the angle formed between the holomorphic complex line and the the curve of complex tangents.

\end{abstract}
\maketitle

% ----------------------------------------------------------------

\section*{I. Introduction}
Complex tangents to an embedding $M^k \hookrightarrow \mathbb{C}^n$ are points $x \in M$ so that the tangent space to $M$ at $x$ (considered as a subspace of the tangent space of $\mathbb{C}^n$) contains a complex line. If $k > n$, all points of $M$ are necessarily complex tangent, by virtue of the dimensions. If $k \leq n$, some (or all) points of $M$ may have strictly real tangent space. If all points of $M$ are real, we say that the embedding is totally real. In general, the dimension of the maximal complex tangent space of $M$ at $x$ is called the dimension of the complex tangent. An embedding is called CR (Cauchy-Riemann) if the complex dimension of all points of $M$ are the same. In [5], M. Gromov used the h-principle to demonstrate that the only spheres admitting totally real embeddings $S^n \hookrightarrow \mathbb{C}^n$ are $S^1$ and $S^3$.
\par\ \par\
In his paper [4], F. Forstneric showed that every closed oriented 3-manifold can be embedded totally real into $\mathbb{C}^3$. In my paper [3], I demonstrate a lack of such flexibility in higher dimensions, particularly we show that the 5-sphere $S^5$ cannot be CR-embedded into $\mathbb{C}^4$. We also show analogous results for different dimensions.
 \par\ \par\
In this work, we will focus on embeddings of closed 3-manifolds $M \hookrightarrow \mathbb{C}^3$. In this situation, all complex tangents are necessarily of dimension $1$ and generically arise along a link. In my paper [2], we demonstrate that for every knot (or link) type in $S^3$ there exists an embedding $S^3 \hookrightarrow \mathbb{C}^3$ which assumes complex tangents exactly along a link  of the prescribed type (with one point being degenerate). In 2018, Kasuya and Takase proved that given any smooth link $L$ in a closed 3-manifold $M$ with trivial fundamental class $[L] \in H_1(M)$, there exists a smooth embedding $M \hookrightarrow \mathbb{C}^3$ whose complex tangents form $\textit{precisely}$ the given link $L$ and are all non-degenerate (see [8]). 
\par\ \par\
We will here generalize the result of Kasuya and Takase to assert that given any smooth orientable 2-plane field over any given  link $L \subset M$ which is null-homologous in $M ([L]=0)$, there exists an embedding $M\ \hookrightarrow \mathbb{C}^3$ whose complex tangents have holomorphic tangent spaces being precisely the specified 2-plane field. In effect, we show that any potential configuration of complex structure to an embedding of a 3-manifold $M$ can actually arise to an actual embedding; the only restriction being that the set of complex tangents $L$ must be null-homologous.
\par\ \par\
We will demonstrate this new result by demonstrating the flexibility to alter the given embedding of Kasuya-Takase to set the holomorphic tangent bundle as we see fit along any given null-homologous link.
\par\
\section*{II. Known Results}
\par\ \par\
The result of Kasuya-Takase in [8] demonstrates that given a closed 3-manifold $M$ and a null-homologous link $L \subset M$, we may construct a smooth embedding $M \hookrightarrow \mathbb{C}^3$ which is complex tangent precisely along $L$. This is in contrast with my earlier result in [1], where I show that any topological type of link may arise as the set of complex tangents to an embedding $S^3 \hookrightarrow \mathbb{C}^3$, with the caveat that the embedding is only $\mathcal{C}^n$ ($n < \infty$) and must have one complex degenerate point. 
\par\ \par\
We now state the result of Kasuya-Takase formally:
\par\ 
\begin{theorem} ([8]): Let $M$ be a closed orientable 3-manifold and $L \subset M$ a closed 1-dimensional submanifold. Then there exists a smooth embedding $f: M \hookrightarrow \mathbb{C}^3$ with complex tangencies exactly along $L$ if and only if $[L] \in H_1(M)$ is trivial.
      \end{theorem}
\par\ \par\
Their method of proof involves using stable maps into $\mathbb{R}^2$ and immersion lifts. In particular, they first demonstrate that given any such null-homologous link $L \subset M$, there exists a liftable stable map $f:M \rightarrow \mathbb{R}^2$ which is singular exactly along the link $L$. Writing the immersion lift as
$(f,g):M \rightarrow \mathbb{R}^4$, we may assume that the immersion lift is an embedding in a tubular neighborhood of the link (perturbing $g$ if necessary).
\par\ \par\
Then writing $f=(f_1,f_2)$ and $g=(f_3,f_4)$, the following smooth immersion is constructed: $G=(f_1,f_2,f_3,f_4,f_1,-f_2):M \hookrightarrow \mathbb{C}^3$ whose complex tangents they show forms precisely the link $L$.
\par\ \par\
With the fact that the immersion is an embedding in a tubular neighborhood of $L$ and using the relative form of the $h$-principle for embeddings, the existence of a smooth embedding $M \hookrightarrow \mathbb{C}^3$ complex tangent exactly along $L$ is assured, and the theorem is proven. (see Gromov in [5] and Kasuya-Takase in [8]) 
\par\ \par\
\section*{III. Specifying Complex Configurations}
\par\ \par\
Let $M$ be a closed 3-manifold and let $L \subset M$ be a link, null-homologous in $M$.  By the work of Kasuya and Takase, there exists a smooth embedding $F:M \hookrightarrow \mathbb{C}^3$ so that the complex tangents of the embedding form exactly $L$. Let $\tau_x \subset T_x M$ denote the holomorphic tangent space at $x \in L$ to this embedding. We wish to show that for any oriented 2-plane field $\eta \in \Gamma(Gr_2 (TM|_L))$ (section of the Grassmann bundle) smoothly varying over $L$ there exists a smooth embedding which is complex tangent exactly on $L$ with holomorphic tangent space being exactly $\eta$. 
\par\ \par\
\begin{theorem}: Let $M$ be a closed orientable 3-manifold and $L \subset M$ be a link with $[L] \in H_1(M)$ trivial and let $\eta$ be a smooth oriented 2-plane field over $L \subset M$. Then there exists a smooth embedding $F: M \hookrightarrow \mathbb{C}^3$ with complex tangencies exactly along $L$ all of which are non-degenerate and with holomorphic tangent space at each $x \in L$ being  $\eta_x$  
      \end{theorem}
\par\ \par\ 
$\underline{\bold{\it{Remark}}}:$
Let $F:M \hookrightarrow \mathbb{C}^3$ be an embedding of the (real) 3-manifold $M$ into complex Euclidean 3-space. Consistent with the previous notation, we consider $\eta_x \subset T_x(M)$ as a real 2-plane in the tangent space for each $x \in L$. By basic differential topology, it is easee that since $F$ is an embedding, we have $M \cong F(M)$ (diffeomorphic) and that correspondingly $T_x (M) \cong T_{F(x)} (F(M))$ as linear subspaces of $\mathbb{C}^3 \cong \mathbb{R}^6$. (see [6])
\par\ \par\
In fact, denoting this induced isomorphism  $dF=F_* : T_x(M) \rightarrow T_{F(x)} F(M)$ as the push forward of the tangent bundle we will show that the holomorphic tangent spaces also coincide under this push-forward map, with an additional assumption we will determine to  be that $dF$ commutes with the complex structure. Note that:
\par\
$$\eta_x = T_x(M) \cap J(T_x(M))$$ (where J is the complex structure on $\mathbb{C}^3$)
\par\ \par\
will form exactly the holomorphic line bundle to $M$ at every point   $x \in L$. We will now demonstrate that the holomorphic tangent bundle on $F(M)$ will be just the push-forward under $F_*$ of the holomorphic tangent bundle assuming that $dF$ commutes with the complex structure,  for a certain class of functions $F$ satisfying certain conditions, and we will have the holomorphic tangent space to $F(M)$ at $F(x)$ written as $\tau_x$ will be given by:
\par\
$$\tau_x = T_{F(x)} ((M)) \cap J(T_{F(x)} (F(M))) = F_*(\eta_x)$$
\par\ \par\
We can in fact show this directly by demonstrating that these two sets (linear spaces) are equal, particularly let $v \in T_x(M) \cap J(T_x(M)) = \eta_x$. Hence:
\par\
$$v \in T_x (M), v \in J(T_x (M)) \implies F_* (v) \in F_*(T_x(M)) \cap F_*(J(T_x(N))$$ 
\par\ \par\
Therefore:
$$f_*(\eta_x) \subset F_*(T_x(M)) \cap F_*(J(T_x(M))$$ 
\par\ \par\
As discussed above, the corresponding tangent spaces are isomorphic (since f is a diffeomorphism) in that $F_*(T_x M) \cong T_{F(x)} (F(M))$  and so the holomorphic tangent space to $F(M)$ at $F(x)$ is given by: $T_{F(x)} (F(M)) \cap J(T_{F(x)}(F(M)) \cong T_x M \cap J(T_x M) = \eta_x$. 
\par\ \par\
Now to directly complete the argument that the sets are equal as above, we need only to show that $T_{F(x)}((F(M)) \cap J(T_{F(x)}((F(M)) \subset F_*(T_x(M)) \cap J(F_* (T_x(M))$. However, to prove this directly we will need the further assumption that each of the matrices $dF_x$ commute with the complex structure $J$. We will prove this commutativity with the complex structure is sufficient to complete the proof of the equality of these sets. See proof of the Claim in page $7$ for the concise proof of this statement.
\par\ \par\
For now, we note that both these linear complex subspaces are necessarily one dimensional as $x$ is a a one dimensional non-degenerate complex tangent and so the holomorphic tangent space $\eta_x$ is a complex line and by the diffeomorphism of the push-forward on the tangent spaces, we see that the holomorphic tangent space to $F(M)$ at $F(x)$ must be a complex line since x is a $1$-dimensional complex tangent to $M$ itself (see [9]). But we saw that $F_* (\eta) \subset T_{F(x)} (F(M)) \cap J(T_{F(x)} (F(M))$ so that the push-forward of the holomorphic tangent space of $M$ at $x$ is a subspace of the holomorphic tangent space to $F(M)$ at $F(x)$ but the complex linear spaces are of the same dimension.
\par\ \par\
Under the necessary assumptions of the commutativity of the differential with $J$, we will demonstrate that the push forward of the holomorphic tangent bundle to $M$ at $x$ will be exactly equal to the holomorphic tangent bundle of the image $F(M)$ at $F(x)$ assuming the conditions needed to make certain that the matrices $f(x)$ commute with the complex structure at every point $x \in L$ (see proof of Claim in page $7$).
\par\ \par\
$\bold{\underline{\it{Proof}}}$: (of Theorem 2)
\par\ \par\
Let $\eta$ be any smoothly varying 2-plane field along $L \subset M$. Consider bases $\{v_x, w_x\} \in T_x(M)$ for $\eta_x$ for each $x \in L$, i.e. a pair of vector fields over $L$ for which $span(v_x, w_x) = \eta_x \subset T_x(M)$ for each $x \in L$ and suppose the consider a complex line field $\tau_x$ to be spanned (over complex numbers) by $\mu_x$. Note then that the real $span(\mu_x, J\mu_x)=\tau_x$. We will now construct matrices $m_x \in GL_6(\mathbb{R}$ so that $m_x \cdot \tau_x = \eta_x$ for each $x \in L$, i.e. matrices which satisfy the equations:
\par\ 
$$m_x (\mu_x) =v_x$$
\par\
$$m_x(J\mu_x) =w_x,  \forall x \in L$$.
\par\ 
This imposes twelve linear conditions on the (real) 36-dimensional space $GL(\mathbb{R}^6)$ for each $x \in L$. We thus obtain 24-dimensional real  space of matrices for each $x \in L$ satisfying these conditions. Denote these spaces by:
\par\
$$\mathcal{S}_x = \{m \in GL_6 (\mathbb{R}) |  m \cdot \tau_x = \eta_x \} \subset GL_6 (\mathbb{R})$$
\par\ \par\
 We will assume without loss of generality that $L = S^1$ is parameterized by the angles:
 $\theta \in [ 0, 2 \pi ]$. 
 \par\ \par\
 Taking $f\theta)$ to be the matrices $m_\theta$ along this curve, we wish to find an embedding $F:M \hookrightarrow \mathbb{C}^3$ so that the derivative for the embedding at each point will be given by the prescribed matrices $m_\theta$, in that $dF_\theta = f(\theta)$ along the curve.
\par\ \par\
Now to find the function $F(\theta)$ we need only integrate $f(\theta)=m_\theta$ along the curve. Note that since $f(0)=f(2\pi)$ as the curve is chosen to be a circle (assumed without loss of generality).
We so define:
\par\ 
$$F(\theta) = \int_0^{\theta} m_\tau d\tau $$
\par\ \par\
Clearly then $F(0)= = \int_0^{0} m_\tau d\tau=0$. The condition that this function is well defined is exactly that $F(2\pi)=G(0)=0$ and as such it is necessary that:
\par \par\ 
$$\int_0^{2\pi} m_\tau d\tau=0$$
\par\ \par\
To ensure the existence of the function $F(\theta)$ that will serve as a well-defined anti-derivative of $f(\theta)$:
\par\
$$dF = m_\theta = f(\theta)$$
\par\ \par\
and have the embedding $F$ with Jacobian along the curve being the given prescribed matrices $m_\theta$ which in turn alter the holomorphic lines in the tangent spaces along the curve. The existence of this embedding $F$ is ensured by our condition that the integral of $f$ along the curve vanishes.
\par\ \par\
Our assumption that:
\par\ \par\
$$\int_0^{2 \pi} m_\tau d\tau = 0$$
\par\ \par\
will also allow us to extend the map $F(\theta)$ accordingly to a local tubular  neighborhood of the link, since the obstruction to such a local extension is exactly the value of this same  integral of $f(\theta)$ about the curve. This corresponds exactly the condition that the $f(\theta)$ may arise as the Jacobian of some embedding $F$ near $L$. This necessary assumption of the vanishing of the integral of $f\theta)$ along the curve is doubly important to ensure both the existence of the well-defined anti-derivative $F(\theta)$ and also for the local extension of $f$ to a tubular neighborhood.
\par\ \par\
This criterion imposes conditions on:
\par\
$$f(\theta) = m_\theta, \theta \in [0, 2 \pi]$$
\par\ \par\
that give rise to $6$ real linear equations on the components of $m_\theta$ that need  be satisfied. The vanishing of this integral of $f(\theta)$ over $[0. 2\pi]$ would then ensure the integrability of $f(\theta)$ which would then serve as the Jacobian of  the embedding $F(\theta)$ along $L$. The value of the integral of $f$  over $[0, 2 \pi]$ also serves as the obstruction to the extension of $f$ to a tubular neighborhood of $L$, which by our assumption is $0$. 
\par\ \par\
(See Hirsch in [7], section 3.1 and Hatcher in [6], section 4.3 as reference to the study of topological extensions of maps and the obstructions there-in)
\par\ \par\
 With the assumption the vanishing of this integral obstruction, we can further extend $F$ to all of $\mathbb{C}^3$ by using a partition of unity beyond the tubular neighborhood of $L$ that is specified by the local extension garnered as above. We can now construct a $\it{smooth}$ function $E:\mathbb{C}^3 \rightarrow \mathbb{C}^3$ so that $dE_x = m_x$, where $m_x \in \mathcal{S}_x$ for each $x \in L$. Since the respective 2-plane fields vary smoothly in $x$ by construction, it is clear that the matrices $m_x$ will also vary smoothly in x. (see Hirsch in [7], section 4.5 as reference for extensions and tubular neighborhoods)
\par\ \par\
We need as above the condition that:
\par\ 
$$ \int_0^{2 \pi} m_\tau d\tau = 0$$
\par\ \par\
As discussed in the Remark on page 3, to ensure that the holomorphic tangent bundle on the image $f(M)$ is just the push-forward of the holomorphic bundle on $M$, we will need that the matrices $m_x$ commute with the complex structure $J$ at each $x \in L$. Note the complex structure $J$ takes the form of a $6x6$-matrix over $\mathbb{R}^6 ( \cong \mathbb{C}^3)$) as follows
\par\
  
$$J= \begin{pmatrix}
0 & -1 & 0 & 0 & 0 & 0  \\
1 & 0 & 0 & 0 & 0 & 0     \\
0 & 0 & 0 & -1 & 0 & 0    \\
0 & 0 & 1 & 0 & 0 & 0     \\
0 & 0 & 0 & 0 & 0 & -1    \\
0 & 0 & 0 & 0 & 1 & 0J
\end{pmatrix}$$
\par\ \par\ \par\
Let $x \in L$ represent a complex tangent as in our previous construction. We have a 24-dimensional space of matrices  $\mathcal{S}_x \subset Gl_6 (\mathbb{R})$ which satisfy the conditions as we gave before in that $m_x$ maps $\tau_x$ to $\eta_x$. We now want to further impose the condition that $m_x$ commutes with the complex structure $J$ for all $x\in L$. Namely we have 18 more real equations on $m_x$ given by the entries of the matrix equation: 
\par\ \par\
$$J \cdot m_x = m_x \cdot J$$
\par\ \par\
More precisely we will have that:
\par\
$$\begin{pmatrix}
0 & -1 & 0 & 0 & 0 & 0  \\
1 & 0 & 0 & 0 & 0 & 0     \\
0 & 0 & 0 & -1 & 0 & 0    \\
0 & 0 & 1 & 0 & 0 & 0     \\
0 & 0 & 0 & 0 & 0 & -1    \\
0 & 0 & 0 & 0 & 1 & 0
\end{pmatrix}
\times
\begin{pmatrix}
a_{11}& a_{12} & a_{13} & a_{14} & a_{15} & a_{16}  \\
a_{21} & a_{22}  & a_{23}  & a_{24} & a_{25} & a_{26}    \\
a_{31}  & a_{32} & a_{33} & a_{34} & a_{35} & a_{36}    \\
a_{41} & a_{42} & a_{43} & a_{44} & a_{45} & a_{46}     \\
a_{51} & a_{52} & a_{53} & a_{54} & a_{55} & a_{56}    \\
a_{61} & a_{62} & a_{63} & a_{64} & a_{65} & a_{66}
\end{pmatrix}$$
$$= \begin{pmatrix}
a_{11}& a_{12} & a_{13} & a_{14} & a_{15} & a_{16}  \\
a_{21} & a_{22}  & a_{23}  & a_{24} & a_{25} & a_{26}    \\
a_{31}  & a_{32} & a_{33} & a_{34} & a_{35} & a_{36}    \\
a_{41} & a_{42} & a_{43} & a_{44} & a_{45} & a_{46}     \\
a_{51} & a_{52} & a_{53} & a_{54} & a_{55} & a_{56}    \\
a_{61} & a_{62} & a_{63} & a_{64} & a_{65} & a_{66}
\end{pmatrix}
\times \begin{pmatrix}
0 & -1 & 0 & 0 & 0 & 0  \\
1 & 0 & 0 & 0 & 0 & 0     \\
0 & 0 & 0 & -1 & 0 & 0    \\
0 & 0 & 1 & 0 & 0 & 0     \\
0 & 0 & 0 & 0 & 0 & -1    \\
0 & 0 & 0 & 0 & 1 & 0
\end{pmatrix}$$
\par\ \par\ \par\ \par\
We can directly see that the only equations that will need to be satisfied are determined by the skew-symmetric behavior:  $a_{ij} = -a_{ji}$ for all $i \neq j, 1 \leq i, j \leq 6$. This will give 15 equations resulting from the choice of two from the six rows, i..e. ${C_2}^6$ combinations of two indices from 6. In addition, we will have three more equations resulting from the necessary equal diagonal elements for each of the "complex pairs":
\par\ \par\
$$a_{11} = a_{22}, a_{33} = a_{44}, a_{55} = a_{66}$$
\par\ \par\
to give a total of 18 equations on the 24-dimensional space $\mathcal{S}_x$ for each complex tangent $x \in \eta$. Let us now write (without loss of generality)  $\mathcal{S}_x$ to be the 6 dimensional space of matrices which satisfy all the earlier conditions, and also commute with the complex structure $J$.
\par\ \par\
We will now construct a new embedding $\tilde{F} = E \circ F: M \rightarrow \mathbb{C}^3$ which will be a composition of smooth functions and hence itself smooth:
\par\ \par\
$\underline{\bold{Claim}}:$  The sets $m_x(T_x M \cap J(T_x M)) = T_x (\tilde{F}(M)) \cap J(T_x (\tilde{F}(M))) $ for every $x \in L$.
\par\ \par\  
First, we note that since $d \tilde{F}|_x = m_x \cdot dF_x$, we have $T_x (\tilde{F}(M)) = m_x T_x M $ for each $x \in L$. By dimensionality arguments, it will suffice to show:
\par\ \par\
$$m_x(T_x M \cap J(T_x M)) \subset m_x T_x M \cap J(m_x T_x M)$$
\par\ \par\
Let $\varphi \in T_x M \cap J(T_x M) = \tau$. Clearly, $\varphi \in T_xM$. So it remains to show that $m_x \varphi \in J(m_x T_x M)$.
\par\ \par\
Since $m_x \in \mathcal{S}_x$, it will commute with the complex structure $J$ by construction, and hence $J(m_x \cdot T_x(M)) = m_x \cdot J(T_xM)$  we see that $m_x(T_x M \cap J(T_x M)) \subset m_x T_x M \cap J(m_x T_x M)$ and by linear algebra dimensionality reasons, the sets must be equal, and the claim is proven.
\par\ \par\
Note that this fact put forth by the claim corresponds with our previous Remark on page $3$ that the $f_*$ push-forward of the holomorphic tangent bundle coincides exactly with the holomorphic bundle of the image $f(M)$.
\par\ \par\
As a direct consequence of this claim, it is clear that the link $L$ will be complex tangent for any such embedding $\widetilde{F}:M \hookrightarrow \mathbb{C}^3$ and the holomorphic tangent space at $x \in L$ is $\eta_x$.
\par\ \par\
The condition as stipulated before  in the previous pages for the "doubly necessary" vanishing of the integral:
\par\
$$ \int_0^{2 \pi} m_\tau d\tau = 0$$
\par\ \par\
leads to $6$ real equations, precisely one for each row dot product with the unit positive tangent to $L$ at the point, and we have now explicitly found the $6$ real equations  that need be satisfied, to be free of any obstructions to extending the solution beyond $L$ and also for the anti-derivative $F(\theta)$ to be well-defined.
\par\ \par\
Together with the 18 equations we found to ensure that the matrix $m_x$ commutes with the complex structure, we now have 24 equations on the 24-dimensional space $\mathcal{S}_x, x \in L$ we require to hold on the matrices, which will now in effect guarantee that they can in fact arise as the derivative (Jacobian) of some embedding as previously discussed.  
\par\ \par
It can be easily verified that these two groups of equations are independent and do not imply (or deny) one another. Therefore, by dimensionality arguments in this linear setting of matrices. there must exist solutions to this system of equations, if fact there will be a unique solution for every $x \in L$.
\par\ \par\
It remains for us to prove that we can find a family $\{m_x \in \mathcal{S}_x\}, x \in L$ so that the resulting embedding $\tilde{F}$ will not have any new complex tangents off of $L$. Note that it will suffice to prove this is true for a tubular neighborhood of $N \subset M$ of $L$ and then proceed using Gromov's h-principle for extensions (see [5]) to extend our solution beyond the tubular neighborhood above to all of $M$. 
\par\ \par\
Define the set: $\mathcal{S} \equiv \cup \mathcal{S}_x  \subset  GL_6 (\mathbb{R})$ over all $x \in  L$. This forms a fiber bundle over $L$ with fibers $\mathcal{S}_x$ over each $x \in L$. Consider further the set of  sections to this bundle  $\bold{\Gamma }=  \{ s: L \rightarrow \mathcal{S}\}$. By our work above, for each section $s$ there exists an automorphism $E:   \mathbb{C}^3  \rightarrow \mathbb{C}^3$ so that the function $\tilde{f}_s=E_s  \circ  f$ has complex tangents along $L$ with holomorphic tangent space exactly the desired complex line field $\eta$. Also, let $N$ be a tubular neighborhood of $L \subset M$.
\par\ \par\
Consider now the "total"  (collection of) Gauss map: $\mathcal{G}:N  \times \bold{\Gamma} \rightarrow Gr_{3, 6}$  that sends $(x, s)$ to the tangent space of the embedding $\tilde{f}_s$ at $x$. We wish to show that $\mathcal{G}$ is transverse to $\mathbb{W} =  \{P  \subset Gr_{3, 6} | P \cap JP \cong \mathbb{C}\}$, the set of planes which contain a complex line.
\par\ \par\
It is readily seen that every direction $v \in \mathbb{R}^6$ can arise in the image of the differential ${dE^s}_x = m_x$ for some $s \in \bold{\Gamma}$, for any given $x \in L$. Hence, we see that the differential $d \mathcal{G}$ is onto and $\mathcal{G}$ is thus necessarily transverse to $\mathbb{W}$. (see [2])  
\par\ \par\
By the parametric transversality theorem, for almost any $s \in \bold{\Gamma}$ , the map $\tilde{f}_s$ will be generic, in that it is either totally real or assumes complex tangents along a link. As $\tilde{f}$ assumes complex tangents along $L$, we can take $N$ to be "sufficiently thin" so that no other knots will be obtained in the transverse intersection. As he intersection of the Gauss map with $\mathbb{W}$ is transverse for such a generic $s$, the complex tangents are all non-degenerate and will form exactly the link $L$. We may then use the totally real h-principle for extensions to construct a smooth embedding $M \hookrightarrow \mathbb{C}^3$ which is (non-degenerate) complex tangent along $L$, and totally real off of $L$, with holomorphic tangent spaces being the arbitrary 2-plane field $\eta$ along $L$, and our theorem is proven (see [5])
\par\ \par\
Note that I used analogous arguments for my proof in [3].
\par\ \par\
$\bold{\it{QED}}$
\par\ \par\
We believe that this result demonstrates the full (topological) flexibility of real embeddings into complex Euclidean space in the dimension 3. It may also be interesting to now ask if the same kind of flexibility will extend to the structure of the holomorphic hulls for real 3-dimensional submanifolds of complex space. 
\par\ \par\

\section*{IV. A Topological Formulation for Analyticity}

The Bishop Invariant assigns to each (non-degenerate) one-dimensional complex tangent $m \in M$  a number $\gamma_m \in [0, \infty]$. Errett Bishop formulated the invariant in the 1950's  (see [1]) for any embedding of a closed real n-manifold $M \hookrightarrow \mathbb{C}^n$ using a certain normal form for $M$ at a one-dimensional complex tangent $m$.  This invariant gives rise to the concept of the $analyticity$ of the complex tangent, and how the manifold may be considered to be (locally) elliptic, hyperbolic, or parabolic near the given complex tangent. It also has implications regarding the holomorphic hull of $M$. The definition is given via Bishop's normal form, which in our case in dimension $3$ is given in a  local holomorphic coordinate system $\{z_1, z_2, z_3 \}$ as follows: (see Webster in [9])
\par\ 
$$M: R=(r_1, r_2, r_3)=0$$
\par\
$$R= \bar{R}$$
\par\
$$dr_1 \wedge dr_2 \wedge dr_3 \neq 0$$
\par\
$$\partial r_1 \wedge \partial r_2 \wedge \partial r_3 = B dz_1 \wedge dz_2 \wedge dz_3$$
\par\ \par\
It is direct to see from [9] that $B$ is a "volume form" type function and is given by:
\par\ \par\
$$B= \frac{\partial(r_1, r_2, r_3)}{\partial(z_1,z_2,z_3)} = \frac{1}{2} \frac{\partial \bar{F}}{\partial z_1} + O(2)$$
\par\ \par\
It is then clear using linear algebra that the set of points $\eta = \{B = 0, R = 0\}$ form exactly the complex tangents of the embedding. For $m \in \eta$, and $H_m = T_m M \bigcap J T_m M$ the holomorphic tangent space, we can write:
\par\ \par\
$$\mathbb{C} \otimes H_m = H'_m \oplus H''_m, H''_m = \overline{H'_m}$$
\par\ \par\
We then need to choose a vector field $X$ which spans $H'_m$. In our situation the $z=z_1$ coordinate represents exactly the holomorphic line $H_m$, and we may choose $X=\frac{\partial}{\partial z}$. 
\par\ \par\
The Bishop Invariant in this situation may then be defined as follows:
\par\
$$\gamma_m = \frac{1}{2} \abs{\frac{B_z}{B_{\bar{z}}}}(m)$$
\par\ \par\
We say that $m \in \eta$ is $\it{degenerate}$ if $B_z = B_{\bar{z}}= 0$ at $m$ and to be $non$-$degenerate$ otherwise. From [7], we find that if $m \in \eta$ is non-degenerate then $\eta$ has a tangent line at $m$.
\par\ \par\
The holomorphic tangent space $H_m$ at a non-degenerate complex tangent $m$ is given exactly by the complex variable $z=z_1$. The tangent space $T_m (M)$ may then be written as the direct sum of the plane $H_m$ and a real line, which we take to form a new real axis which we may designate as the $x$-axis, and now the tangent space $T_m(M) = \{(z, x) | z \in \mathbb{C}, x \in \mathbb{R} \}$.
\par\ \par\
Consider now the normal bundle $S$ of $\eta \subset M$ and let $S^*$ denote its conormal bundle. It is direct to see that $S^*$ will be spanned by the co-frame $dB$ along $\eta$.
\par\ \par\
Consider then the composite function:
\par\
$$\varphi:H \hookrightarrow T(M)|_\eta \rightarrow  (T(M)|_\eta)/T(\eta) \equiv S$$.
\par\
In reference to Webster in his work in [9], we say that $m \in \eta$ is $\it{elliptic}$ if $\varphi_m$ is orientation-reversing. This will be true if and only if $\gamma_m < \frac{1}{2}$.
\par\
We may also say $m$ is $\it{hyperbolic}$ if $\varphi_m$ is orientation-preserving. This will be true if and only if $\gamma_m > \frac{1}{2}$.
\par\
We may say $m$ is $\it{parabolic}$ if $\varphi_m$ is singular of rank 1. This is equivalent to the situation where $\gamma_m = \frac{1}{2}$.
\par\ \par\
It is then clear that $m$ is parabolic precisely when $T_x (\eta) \subset H_x$. It is also evident that if the tangent line $T(\eta)$ is oriented to be "above" the holomorphic space $H_m$ in agreement with the orientation of the positive normal of $H_m$ then $m$ will be a hyperbolic complex tangent. Analogously, if $T(\eta)$ is oriented "below" and in disagreement with the orientation of $H_m$, then $m$ will be an elliptic complex tangent.
\par\ \par\
If $T(\eta) \perp H_m$, we say that $m$ is an $\it{exceptional}$ complex tangent. In accordance with the definition of the Bishop invariant, the condition that $\gamma=0$ will be the extreme case for ellipticity and will satisfy that $T(\eta) = {\overline{H}_m}^{\perp}$ with opposite orientation to that induced by $H_m$. Analogously, $\gamma = \infty$ will correspond to the extreme hyperbolic case  $T(\eta) = (H_m)^{\perp}$ in agreement with orientation induced by $H_m$. 
\par\ \par\
From our previous work in Section 3, we are able to construct an embedding with the set of complex tangents being along any curve $\gamma \subset M$ with the holomorphic planes along any 2-plane field $\xi \subset Gr_2(M)|_\eta$ and as such we can control the flow of the tangents of $\eta$ along the planes $\xi$ and by the above we can form the embedding to have any "Bishop structure" along the complex tangents as we so choose.
\par\ \par\ \par\
We may then consider the angle $\theta_m$ formed by the unit normal $\rho_m$ to $H_m \subset T_m (M)$ (using the natural orientation of induced by the complex structure) and $t_m \in T_m(M)$ representing the unit positive tangent to the curve of complex tangents $\eta$ at $m$. We then summarize our work above in the classification of the complex tangents of the embedding in terms the angles $\theta_m = \angle (t_m, \rho _m)$ as follows:
\par\ \par\   \par\ \par\ \par\ \par\ \par\ \par\   \par\ \par\ \par\ \par\ \par\ \par\ \par\  

 \begin{theorem}: Let $M$ be a closed orientable 3-manifold and $L \subset M$ a homologically trivial closed 1-dimensional submanifold (fundamental class [L]=0), and let  $\eta$ be a smooth oriented 2-plane field over $L \subset M$.  Then there exists a smooth embedding $M \hookrightarrow \mathbb{C}^3$ non-degenerate complex tangent exactly along $L$ with holomorphic tangent space $\eta_x$ for each $x \in L$. Furthermore, for each $m \in L$ the analyticity of $m$ may be (pre-) determined by the angle $\theta_m = \angle (t_m, \rho _m)$ between the complex line and the curve of complex tangents as follows:
\par\ \par\ \par\
1.	If $\theta_m \in [0, \pi /2)$ then $m$ is a hyperbolic complex tangent.
\par\ \par\
2.	If $\theta_m \in (\pi /2, \pi]$ then $m$ is an elliptic complex tangent.
\par\ \par\
3.	If $\theta_m = \pi /2$ then $m$ is a parabolic complex tangent.
\par\ \par\
If $\theta_m \in \{0, \pi\}$ then $m$ is an exceptional complex tangent, of elliptic type if $\theta_m = \pi$ and of hyperbolic type if $\theta_m = 0$.
      \end{theorem}

 \par\ \par\ \par\
We also note that it could be of interest to find a new formula for the Bishop invariant in terms of the angles $\theta_m$ among the set of complex tangents $m \in N$. We will leave  this work for another paper as it would need an analytic sophistication that we do not make use of in our current work, which has been in the spirit of complex differential topology in nature.
\par\ \par\ \par\ \par\

\end{document}